\documentclass[final]{siamonline1116}

\usepackage{lipsum}
\usepackage{amsfonts}
\usepackage{graphicx}
\usepackage{epstopdf}
\usepackage{algorithmic}
\usepackage{subcaption}
\ifpdf
  \DeclareGraphicsExtensions{.eps,.pdf,.png,.jpg}
\else
  \DeclareGraphicsExtensions{.eps}
\fi

\numberwithin{theorem}{section}

\newcommand{\TheTitle}{Effects of periodic forcing on a Paleoclimate delay model} 
\newcommand{\TheAuthors}{C. Quinn, J. Sieber, and A. S. v. d. Heydt}

\headers{\TheTitle}{\TheAuthors}

\title{{\TheTitle}\thanks{Submitted to editors 11/1/19.
\funding{C.Q. and J.S have received funding from the European Union's
Horizon 2020 research and innovation programme under the Marie
Sklodowska-Curie grant agreement No 643073. J.S. gratefully
acknowledges the financial support of the EPSRC via grants
EP/N023544/1 and EP/N014391/1. AvdH gratefully acknowledges financial 
support from the EPSRC funded Past Earth Network (grant number EP/M008363/1) 
and ReCoVER (grant number EP/M008495/1) for an extended research 
visit to the University of Exeter in summer 2017. }}}

\author{
  Courtney Quinn\thanks{College of Engineering, Mathematics, and Physical Sciences, University of Exeter, Exeter EX4 4QE, United Kingdom
    (\email{c.quinn2@exeter.ac.uk}, \email{J.Sieber@exeter.ac.uk}).} 
  \and
  Jan Sieber\footnotemark[2]
  \and
  Anna von der Heydt\thanks{Institute for Marine and Atmospheric Research,
Department of Physics \& Center for Complex Systems Studies,
Utrecht University, Princetonplein 5, 3584 CC
Utrecht, The Netherlands
    (\email{A.S.vonderHeydt@uu.nl}).}
}

\usepackage{amsopn}

\DeclareMathOperator{\leb}{Leb}
\DeclareMathOperator{\rg}{rg}
\DeclareMathOperator{\dom}{dom}

\DeclareMathOperator{\mae}{MAE}
\newcommand{\R}{\mathbb{R}}
\newcommand{\Z}{\mathbb{Z}}
\newcommand{\Lint}{\mathbb{L}}
\graphicspath{{figures_periodic/}}
\definecolor{dgreen}{rgb}{0,0.5,0}
\definecolor{dblue}{rgb}{0,0,0.5}
\definecolor{dred}{rgb}{0.5,0,0}

\ifpdf
\hypersetup{
  pdftitle={\TheTitle},
  pdfauthor={\TheAuthors}
}
\fi

\newcommand{\ree}{$r_{\mathrm{ee}}$}
\newcommand{\res}{$r_{\mathrm{es}}$}
\newcommand{\re}{$r_{\mathrm{e}}$}
\newcommand{\rel}{$r_{\mathrm{el}}$}
\newcommand{\rl}{$r_{\mathrm{l}}$}

\begin{document}

\maketitle

\begin{abstract}
  We present a study of a delay differential equation (DDE) model for
  the Mid-Pleistocene Transition (MPT). We investigate the behavior of
  the model when subjected to periodic forcing.  The unforced model
  has a bistable region consisting of a stable equilibrium along with
  a large amplitude stable periodic orbit.  We study how forcing 
  affects solutions in this region. Forcing based on astronomical data 
  causes a sudden transition in time and under increase of the forcing 
  amplitude, moving the model response from a non-MPT
  regime to an MPT regime.  Similar transition behavior is found for
  periodic forcing.  A bifurcation analysis shows that the transition
  is not due to a bifurcation but instead to a shifting basin of
  attraction. While determining the basin boundary we demonstrate
    how one can accurately compute the intersection of a stable
    manifold of a saddle with a slow manifold in a DDE by embedding
    the algorithm for planar maps proposed by England \emph{et al.}
    (SIADS 2004(3)) into the equation-free framework by Kevrekidis
    \emph{et al.} (Rev. Phys. Chem. 2009 (60)).
\end{abstract}



\section{Introduction}

Low-dimensional conceptual models are often used in climate modelling to understand basic interactions between specific climate variables \cite{dijkstra2013}.  These are particularly useful when studying long time-scale dynamics and computing power precludes resolving all of the many temporal and spatial scales \cite{kaper2013}.  For this reason, conceptual models are essential in studying past climates of the earth and their long-term variability.

\begin{figure}
        \centering
        \includegraphics[width=0.75\textwidth]{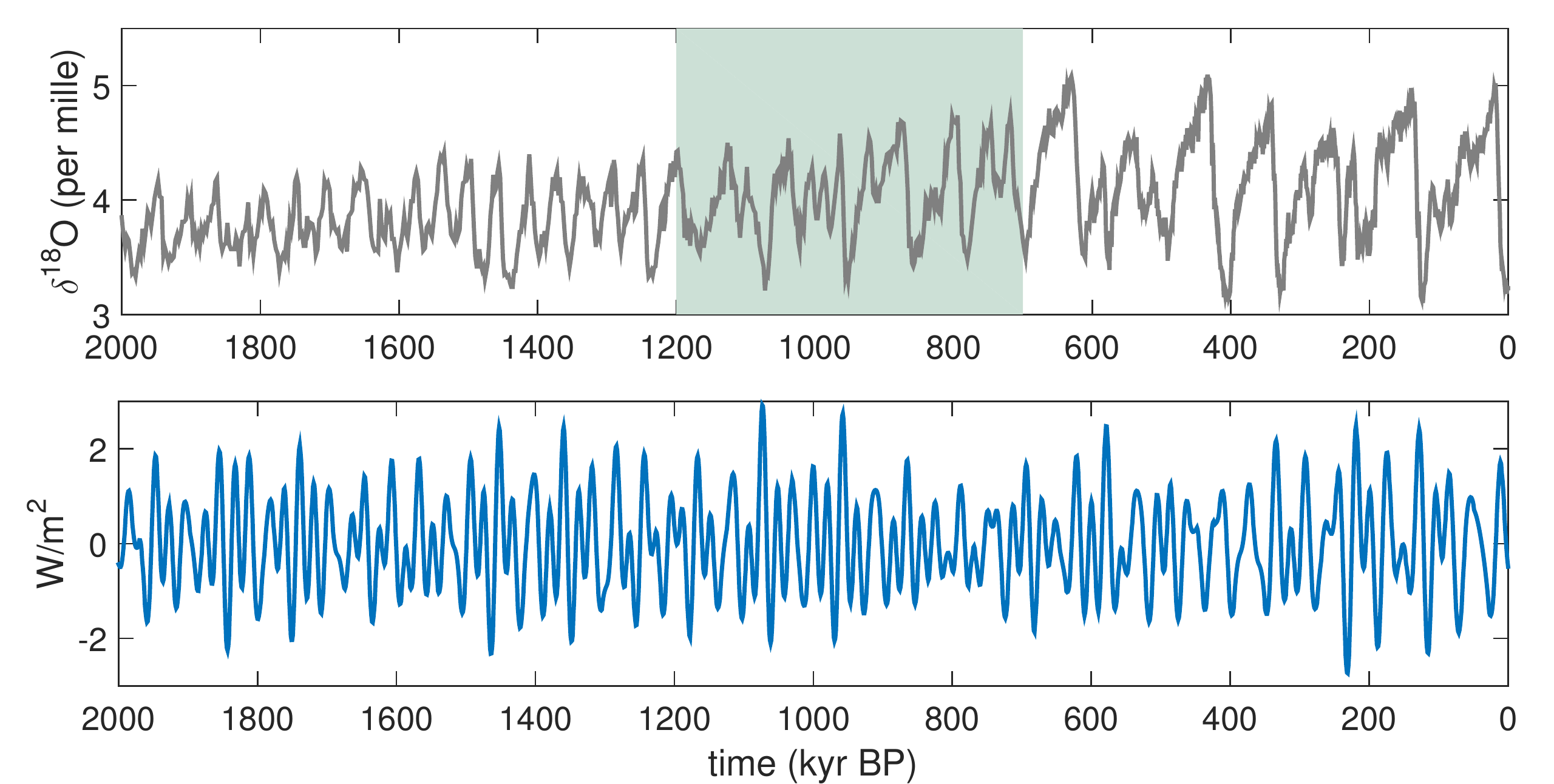}\includegraphics[width=0.24\textwidth]{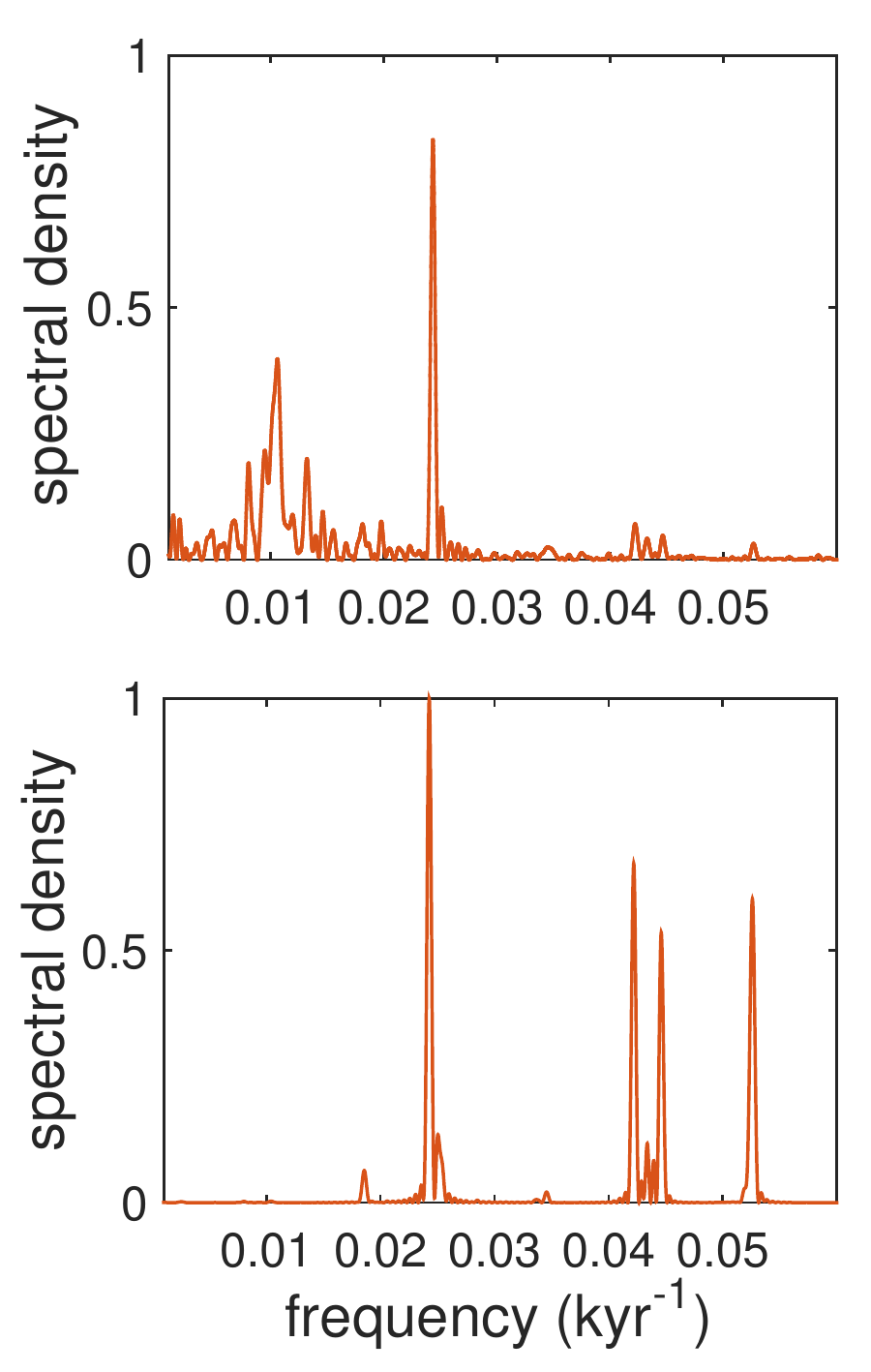}
        \caption{\label{fig:forcing} (Top) Proxy record of global ice cover for the last 2 Myr \cite{lisiecki2005} and its spectrum.  (Bottom) Normalized integrated July insolation $F_I(t)$ at 65$^{\circ}$N, adapted from \cite{huybers2006}, and its spectrum.}
\end{figure}

The Pleistocene, which lasted from approximately 2.6 Myr to 11.7 kyr before present, is a period in climate history that particularly benefits from this type of modelling.  The essential variability can be captured through dynamical systems with a only a few coupled variables (see review \cite{crucifix2012} for an extensive collection of examples).  During this time, the earth experienced many oscillations between periods of colder temperatures and increased global ice cover (\emph{glacials}) and periods of warmer temperatures with less global ice cover (\emph{interglacials}) \cite{imbrie1986}.  The dominant periodicity of the oscillations also changed from approximately 41 kyr in the beginning of the Pleistocene to roughly 100 kyr towards the end of the Pleistocene, together with an increase in amplitude and degree of asymmetry in the oscillations \cite{hays1976,maasch1988}.  This shift in dynamics is known as the Mid-Pleistocene Transition (MPT), and the exact timing of it is believed to be sometime between 1200 and 700 kyr BP \cite{maasch1988,paillard2004,crucifix2012,dijkstra2013,engler2017}.  The oscillations and the MPT can be observed through proxy records as shown in \cref{fig:forcing}.

In a recent study \cite{quinn2017} we revisited a model of the Pleistocene introduced by Saltzman and Maasch \cite{saltzman1988}.  This model was used to replicate the main dynamics of the ice ages, which involved perturbations in glocal ice mass, atmospheric CO$_{2}$, and global ocean circulation. In the original study, the authors of \cite{saltzman1988} showed this particular model was able to reproduce the MPT through a slow parameter shift.  In our study \cite{quinn2017} we reduced this model to a scalar delay differential equation (DDE) for global ice mass.  Through analytical and numerical investigations we showed that the models were qualitatively equivalent.  Our main focus of the study was on the bistable region found in both the original ODE model and the DDE model, a region which was not explored in the original analyses \cite{saltzman1988,maasch1990}.  We observed the behavior within this region when the model was subjected to external forcing, namely solar radiation, and were able to reproduce a MPT-like transition without any change in parameters.

Very early studies have shown a relationship between solar radiation and the glacial cycles \cite{milankovitch1941,budyko1969,sellers1969}. The solar radiation, otherwise known as astronomical insolation, is a quasiperiodic forcing with dominant frequencies around 0.0243, 0.0434, and 0.0526 kyr$^{-1}$, corresponding to periodicities of 41, 23, and 19 kyr respectively (see \cref{fig:forcing}) \cite{berger1978,huybers2006early}.  The most prominent signal in the forcing is at 41 kyr, and this has been argued to be the driving force of the climate fluctuations in the beginning of the Pleistocene \cite{milankovitch1941}, while it remains still unclear where the  longer-term and large-amplitude fluctuations of the late Pleistocene derive from and why a transition from the 41 kyr to the 100 kyr dominant periodicity occurred \cite{Paillard:2015is}. It is natural then to consider the effects of the 41 kyr signal alone on our model.  In this paper we perform a systematic study of how the model responds to periodic forcing with a period of 41 kyr.

The paper is organized as follows. We summarize the numerical observations and the bifurcation analysis of the unforced system from our recent study \cite{quinn2017} in \cref{sec:backgd}. \Cref{sec:main} studies the effects of periodic forcing, showing that one type of transition observed in \cite{quinn2017} can be attributed to a shift in basin boundaries in the context of periodic forcing. We are able to track the precise basin boundary by applying algorithms developed for stable manifolds of two-dimensional maps to the forced DDE. \Cref{sec:conclusions} compares our results to other hypotheses, also suggesting how the theory for quasiperiodically forced systems could be applied to extend our results.

\section{Background}
\label{sec:backgd}

Our delay equation model of global ice mass perturbations as described in \cite{quinn2017} is as follows:
\begin{equation} \label{eq:DDE}
\dot{X}(t) = -pX(t-\tau)+rX(t)-sX(t-\tau)^{2}-X(t-\tau)^{2}X(t).
\end{equation}
Here, $X$ represents the ice mass anomaly from a background state
where $X(t)$ and $X(t-\tau)$ are taken at present and $\tau$ years in
the past, respectively.  The delay $\tau$ is the timescale of the
feedback processes associated with ice accumulation and decay, and
carbon storage and transport in the deep ocean.  This is the parameter
we are most interested in analysing.  The other parameters $p$, $r$,
and $s$ will be kept constant at $p=0.95$, $r=0.8$, and $s=0.8$ in
accordance with \cite{saltzman1988} for all numerical studies. The
  unit for time $t$ and delay $\tau$ is $10$ kyr throughout the paper,
  unless the time unit is explicitly specified (then the time axis is
  usually kyr BP, as in \cref{fig:forcing}).

This is not the first time a delay model has been suggested for the dynamics of the Pleistocene.  Bhattacharya \emph{et al.} \cite{bhattacharya1982} explored an energy balance model which incorporates a delay related to feedback effects from the reflectivity of the earth's surface.  Additionally, Ghil \emph{et al.} \cite{ghil1987} proposed a Boolean delay model for global temperature, northern hemisphere ice volume, and deep-ocean circulation, with delays corresponding to ice sheet expansion, ice accumulation, and overturning time of the deep ocean.  The three delay effects discussed in \cite{ghil1987} are captured within our delayed feedback model.

\subsection{Internal Dynamics}
\begin{figure}
        \centering
        \includegraphics[width=0.74\textwidth]{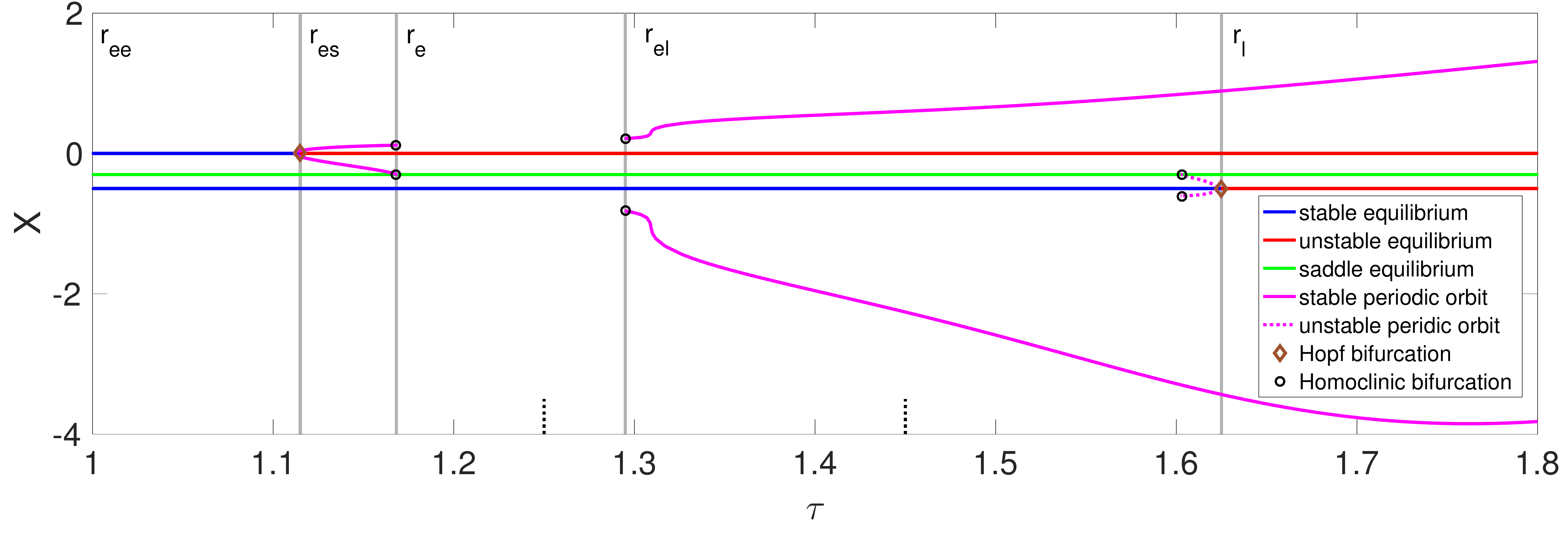} \includegraphics[width=0.25\textwidth]{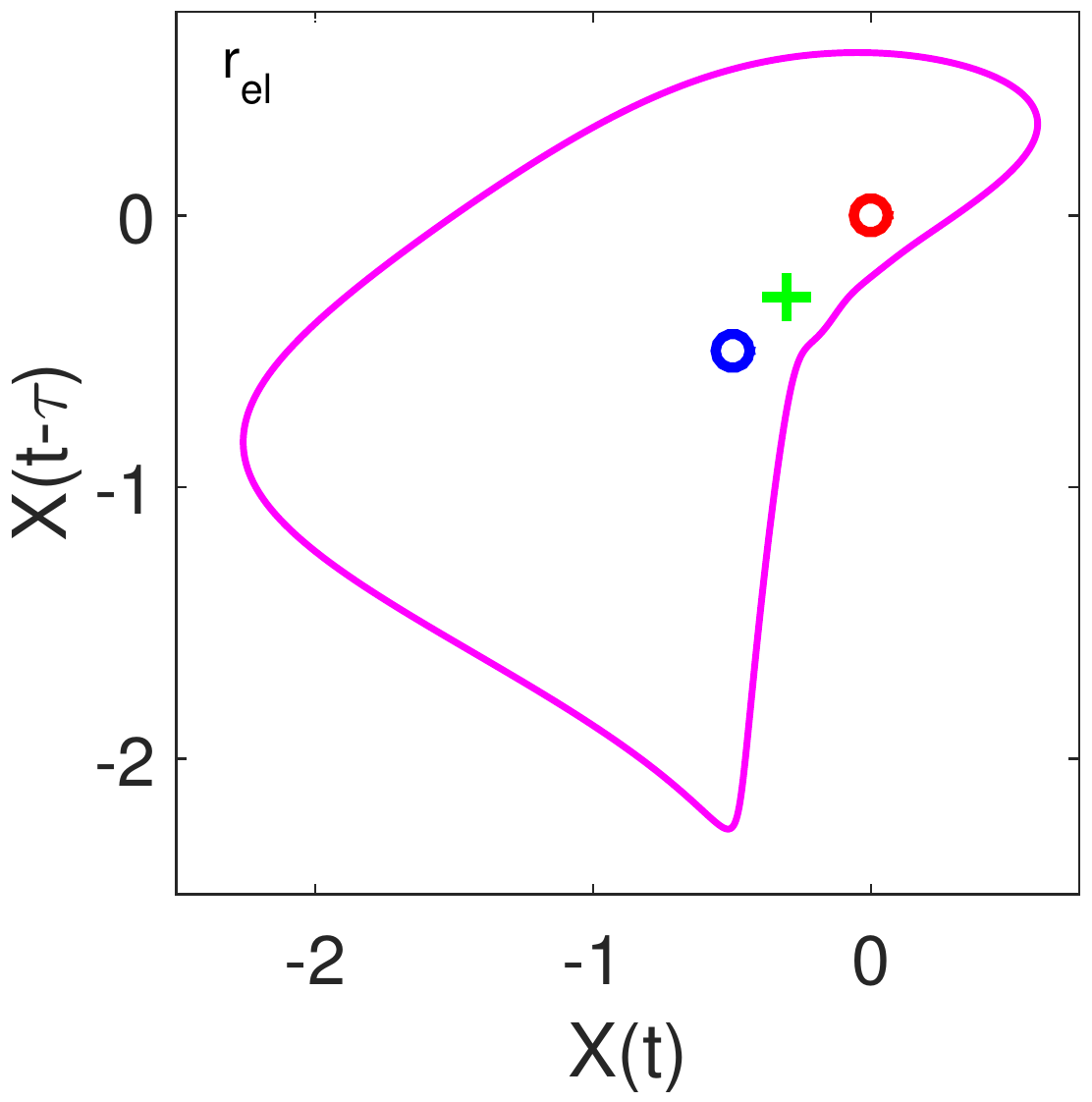}
        \caption{\label{fig:taubif_DDE} (Left) Bifurcation diagram of the DDE model \eqref{eq:DDE} for delay parameter $\tau$.  The dotted black lines indicate values of $\tau$ used in forcing studies: $\tau_{\mathrm{ref}}=1.25$ and $\tau_{\mathrm{bist}}=1.45$. Figure adapted from \cite{quinn2017}. (Right) Phase portrait of the bistable region r$_{\mathrm{el}}$, $\tau_{\mathrm{bist}}=1.45$. The circles are the stable (blue) and unstable (red) equilibria, the green cross is the saddle equilibrium, and the pink curve shows the periodic orbit. Other parameters: $p=0.95$, $r=s=0.8$.}
\end{figure}

We conducted a bifurcation analysis of the model for realistic values of the delay, $\tau\in(1,2)$.  There were five distinct regions with respect to global stability (see \cref{fig:taubif_DDE}).  They are as follows:  
\begin{itemize}
\item {[\ree]} two stable equilibria,
\item {[\res]} one stable
  equilibrium and one stable small-amplitude periodic orbit,
\item {[\re]} one stable equilibrium,
\item {[\rel]} one stable equilibrium
  and one stable large-amplitude periodic orbit, and
\item {[\rl]} one stable large-amplitude periodic orbit.
\end{itemize}
The bistable region, \rel, for $\tau\in[1.295,1.625]$ with a stable equilibruim and a large amplitude stable periodic orbit was previously not explored. Within this region, if the model is subjected to external forcing, transitions are possible between the two stable states without any change in parameters.

We also make a note about the dimensionality of the system. Although DDEs are infinite-dimensional, the phase portrait in \cref{fig:taubif_DDE} gives initial evidence that the dynamics are confined to a two-dimensional slow manifold.  Engler \emph{et al.} \cite{engler2017} derived a two-dimensional slow manifold for the original model \cite{saltzman1988} through considering the deep ocean timescale as instantaneous ($\tau=1$).  Here we consider the case where the deep ocean timescale is not instantaneous ($\tau>1$).  The dimensionality will be investigated in more detail in Section~\ref{sec:slow}.

\subsection{Astronomical Forcing leading to sudden transition}
\label{sec:astro}
Quinn \emph{et al}.~\cite{quinn2017} studied the model's response in
the bistable region \rel{} when subjected to astronomical forcing using
simulations.  They included the forcing as an additive term (with
negative amplitude as insolation reduces ice mass),
\begin{equation}\label{eq:DDE_forced}
\dot{X}(t) = -pX(t-\tau)+rX(t)-sX(t-\tau)^{2}-X(t-\tau)^{2}X(t)-u F_I(t). 
\end{equation}
The term $F_I(t)$ is the forcing signal, shown in \cref{fig:forcing},
which is a time series of integrated summer insolation at
65$^{\circ}$N computed by Huybers \cite{huybers2006} based on the
model in \cite{huybers2006early}.  Details, how the forcing curve in
\cref{fig:forcing} was obtained from publicly available data, are
given in \cite{quinn2017}.  The data was obtained through
  numerical approximation of changes in the long-term planetary motion
  based on the theory for long-term variation of daily insolation
    by Berger \cite{berger1978}. As discussed by Huybers
    \cite{huybers2006early}, the forcing is dominated by a
    quasi-periodic superposition of approximately periodic variations
    of precession, obliquity and eccentricity. This signal, shown in \cref{fig:forcing} (bottom
  panel), can be approximated by a quasiperiodic function of 35
  frequencies \cite{crucifix2012}, with dominant frequencies corresponding to periods of
around 41, 23, and 19 kyr (compare to the lower right panel of
\cref{fig:forcing}).  The parameter $u$ represents the forcing
amplitude of which the realistic value is uncertain.

\begin{figure}
\centering
\begin{subfigure}{0.33\textwidth}
        \includegraphics[width=1.1\textwidth]{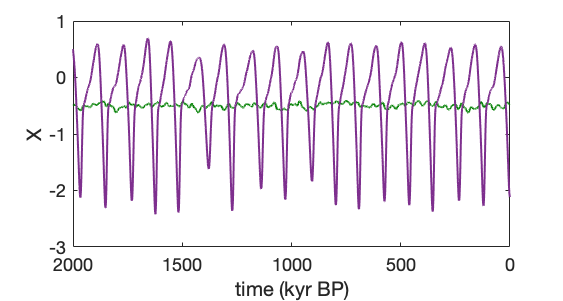}
        \caption{}
        \label{fig:ex_bist_new}
\end{subfigure}
\begin{subfigure}{0.32\textwidth}
        \includegraphics[width=1.1\textwidth]{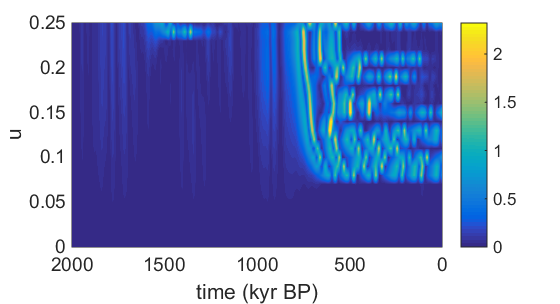}
        \caption{}
        \label{fig:heatmap_tau145}
\end{subfigure}
\begin{subfigure}{0.33\textwidth}
        \includegraphics[width=1.1\textwidth]{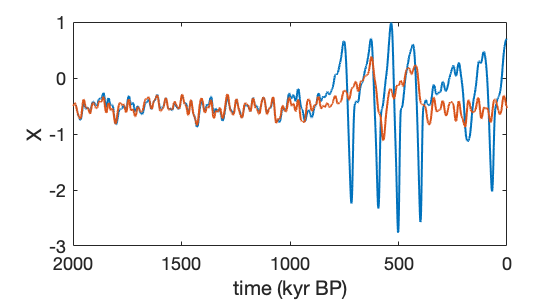}
        \caption{}
        \label{fig:exbist_tau145}
\end{subfigure}
\caption{(adapted from \cite{quinn2017}) \emph{Small-amplitude} and \emph{large-amplitude} responses
  in bistable region \rel{} when subjected to astronomical
  forcing. (a) Example trajectories of small-amplitude (green)
    and large-amplitude response (purple) for
    $\tau=\tau_{\mathrm{bist}}$ and $u=0.05$.  (b) Distance
    between trajectory for delay $\tau_\mathrm{bist}=1.45$ and
    reference trajectory $\tau_{\mathrm{ref}}=1.25$ at given forcing
    amplitude $u\in[0,0.25]$.  Averages taken over window length of
  size $\tau_{\mathrm{bist}}$.  (c) Example of transitioning
    trajectory in bistable region for $\tau_{\mathrm{bist}} = 1.45$
    and $u=0.15$ (blue) and its reference trajectory at
  $\tau_{\mathrm{ref}} = 1.25$ and $u=0.15$ (red).  Initial histories
  $X(s)=-0.5$ for $s\in [2000+10\tau ,2000]$ kyr BP for panels
    (b), (c) and green response in panel (a), $X(s)=+0.5$ for purple
    response in (a); parameters $p=0.95$, $r=s=0.8$ for all
    panels.}
\label{fig:exbist}
\end{figure}
\cref{fig:exbist} summarizes the most important observation of
  Quinn \emph{et al}.~\cite{quinn2017}, which motivates our
  investigation in the following sections.  For values of $\tau$ in
the bistable region \rel{}, there are two possible responses for sufficiently
small $u$, both illustrated in \cref{fig:ex_bist_new}: a
\emph{small-amplitude} (green in \cref{fig:ex_bist_new}) and a
\emph{large-amplitude} (purple) response.  Each response is a
perturbation of an attractor of the unforced system, namely the
equilibrium and the large-amplitude periodic orbit, which persist for
small $u$. The green time profile has been computed starting from a
  constant history of $X(s)=-0.5$ (the autonomous stable equilbrium)
  for $s\in [2000+10\tau ,2000]$ kyr BP (recall that the time unit
    for $\tau$ was $10$ kyr) and shows only the small-amplitude
  response. The purple time profile is computed starting from a constant
  history of $X(s)=+0.5$ for $s\in [2000+10\tau ,2000]$ kyr BP
  and shows the large-amplitude response.
As both responses are perturbations of attractors they will not
  change when we perturb the initial conditions slightly (apart from a
  short transient) \cite{quinn2017}. This persistence breaks down as
  one increases the forcing amplitude
  $u$.  
  Numerical experiments by Quinn \emph{et al.}~\cite{quinn2017}
  discovered that this breakdown leads to a noticeable transition
  between these two reponses. \Cref{fig:heatmap_tau145} summarizes
  the transition for forcing amplitudes $u\in[0,0.25]$ and fixed delay
  $\tau_\mathrm{bist}=1.45$ in the bistable region \rel{} of
  \cref{fig:taubif_DDE}:
  \begin{itemize}
  \item the transition occurs prominently with respect to change of forcing strength $u$ near $u\approx 0.08$,
  \item it occurs with respect to time, consistently between 700
      and 750 kyr BP. This agrees with the timing of the MPT according
    to palaeorecords (see upper panel of \cref{fig:forcing}).
  \item Not visible in \cref{fig:heatmap_tau145}, but discussed in
      \cite{quinn2017}, the transition occurs independent of small perturbations to the
      initial condition.
\end{itemize}
The blue time
  profile in \cref{fig:exbist_tau145} is one example of a solution
  that displays the transition in time (at $u=0.15$).  It
shows a small-amplitude response  up until about
750 kyr BP, when it then transitions to the large-amplitude
response.  For comparison
\Cref{fig:exbist_tau145} also shows the red time profile, which is the
response for identical $u$, but for delay $\tau_\mathrm{ref}=1.25$
outside of the bistable region, where only the small-amplitude
response exists. The color coding in the overview
  \cref{fig:heatmap_tau145} is defined by the distance between the
  responses at $\tau=\tau_\mathrm{bist}=1.45$ and
  reference trajectories $\tau=\tau_\mathrm{ref}=1.25$ with identical forcing strength $u\in[0,0.25]$ and
  identical initial histories (at $X(s)=-0.5$). So, for example, at
  $u=0.15$ the color in \cref{fig:heatmap_tau145} is defined by
  the distance between blue and red profile in
  \cref{fig:exbist_tau145}.

The supplementary material includes a video of trajectories for
$\tau=1.45$ (in blue) and $\tau=1.25$ (reference, red) as they change
with increasing $u$, where more examples of MPT-like transitions can
be seen.  The fact that the system favors this time period to
transition is beyond the scope of this paper, but we will
  describe a possible mechanism for a simplified forcing scenario in
  \cref{sec:conclusions}.

In this study we will focus on the first feature of the
  transition, the transition with respect to forcing strength
$u$. As we will show, the transition in forcing amplitude $u$
  occurs already when the model is subject to periodic forcing. Studying the periodic case first will help us predict
the range of forcing amplitudes that allows a temporal transition
similar to the MPT (possibly timed by a distinct property of the
astronomical forcing around the transition). Since 41 kyr is
  the prominent periodicity in the astronomical forcing, we expect
  that the system shows a similar transition to a large amplitude
  response when varying the forcing amplitude $u$ of a periodic
  forcing with period 41 kyr.

\section{Periodic Forcing}
\label{sec:main}

We are interested in the behavior of this model when a sinusoidal forcing with period 41 kyr is included.  This corresponds to the most prominent frequency found in orbital forcing - the obliquity variations, \textit{i.e.} the changes in the angle between the rotational and orbital axes. Thus, we choose
\begin{equation} \label{eq:forcing_periodic}
F_P(t) = \sin(2\pi t/T).
\end{equation}
The forcing period is  $T= 4.1$, corresponding to $41$\,kyr.

\begin{figure}
\centering
\begin{subfigure}{0.48\textwidth}
\centering
\includegraphics[width=\textwidth]{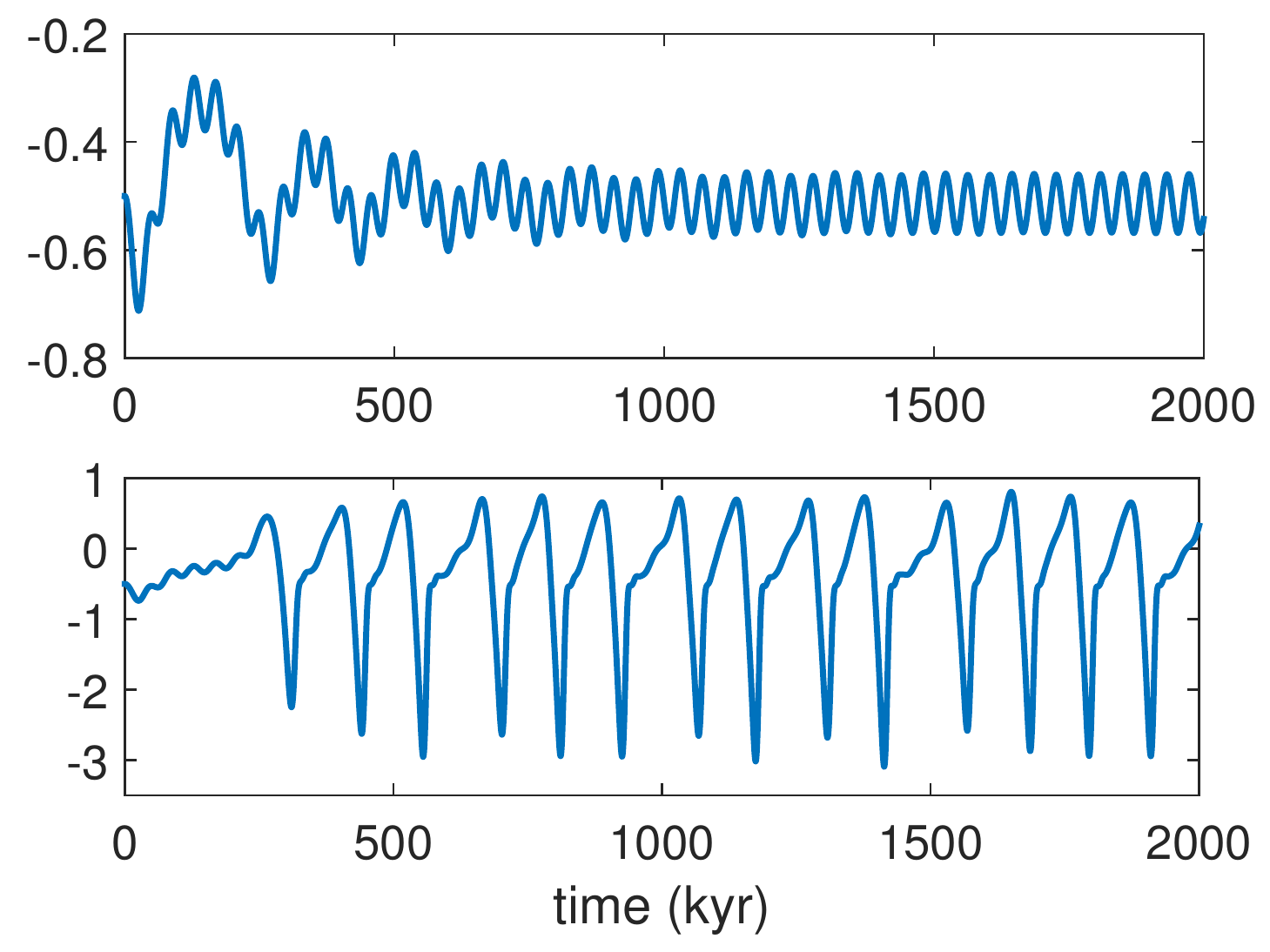}
\caption{Top: $u=0.08$, Bottom: $u=0.09$}
\label{fig:ex_traj} 
\end{subfigure}
\begin{subfigure}{0.48\textwidth}
  \centering
        \includegraphics[width=\textwidth]{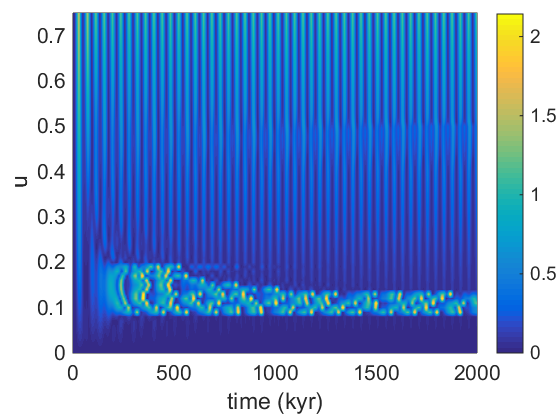}
        \caption{Distance from unforced stable equilibrium}
        \label{fig:heatmap}
\end{subfigure}
\caption{(a) Example trajectories of \eqref{eq:periodic} for two
  qualitatively different long-time regimes. Top: small-amplitude
  response (a small-amplitude periodic orbit), bottom: large-amplitude
  response (longer time series suggest that it is chaotic). (b)
  Distance of solution $X$ from $X_0:t\mapsto-0.5$ (a stable equilibrium
  of the unforced system) for varying forcing amplitudes $u$. Other
  parameters: $\tau=1.55$, $T=4.1$, $p=0.95$, $r=s=0.8$, $\phi=0$;
  initial history $X_0(s)=-0.5$ for $s\in [-\tau,0]$.}
\end{figure}

This forcing is included in the same way as the astronomical forcing $F_I$,
\begin{equation}\label{eq:periodic}
\dot{X}(t) = -pX(t-\tau)+rX(t)-sX(t-\tau)^{2}-X(t-\tau)^{2}X(t)-uF_P(t)\mbox{.} 
\end{equation}
DDE~\eqref{eq:periodic} is a dynamical system with the phase space
$U=C([-\tau,0];\R)$, where $C([-\tau,0];\R)$ is the space of continuous functions on the interval $[-\tau,0]$ with the maximum norm $\|X\|_0=\max\{|X(t)|:t\in[-\tau,0]\}$. At any given time $t\geq0$, the state is $X_t:[-\tau,0]\ni s\mapsto X(t+s)\in\R$.
For sufficiently small values of $u$ and $\tau\in r_\mathrm{el}$,
there exist
\begin{itemize}
\item a stable small amplitude periodic orbit (with period $T$), which is
  a perturbation of the stable equilibrium at $u=0$, and
\item a stable quasiperiodic large amplitude solution, which is a
  perturbation from the large amplitude periodic orbit at $u=0$.
\end{itemize}
These two attractors will persist for a range of $u$ and we will refer
to them as the \emph{small-amplitude response} and the
\emph{large-amplitude response}, as we did in the case of astronomical
forcing. Both types of stable long-time regimes are shown in
\cref{fig:ex_traj} including a transient. We observed in simulations
that the large-amplitude response changes from quasiperiodic to
chaotic as $u$ increases.  Large-amplitude chaotic responses have been
observed prevsiously in conceptual ice age models subject to periodic
forcing in the literature.  Ashwin \emph{et al.} \cite{ashwin2018}
find significant regions of chaotic responses for the van der
Pol-Duffing oscillator, the Saltzman and Maasch 1991 model
\cite{saltzman1991}, and the Paillard and Parrenin 2004 model
\cite{paillard2004}.  The chaos exists both for simple periodic
forcing defined by \cref{eq:forcing_periodic} and more complex
quasiperiodic forcings.  In contrast to our scenario, in
\cite{ashwin2018} all of the models were considered in parameter
regions where the unforced dynamics has a single large-amplitude
stable periodic orbit. Our simulations suggest that large-amplitude
chaotic solutions are also present in a periodically forced bistable
regime.
 
The heat map in \cref{fig:heatmap} shows the model response over a
larger range of forcing amplitudes $u$.  For \cref{fig:heatmap} we
keep the delay constant at $\tau=1.55$ and increase $u$ from $0$ to
$0.75$.  All trajectories start from the constant initial history
$X_0:[-\tau,0]\ni s\mapsto-0.5\in\R$ corresponding to the stable
equilibrium of the unforced system.  We then compute the distance of
$X_t:[-\tau,0]\ni s\mapsto X(t+s)\in\R$ to $X_0$, using the mean
absolute error ($\mae$),
$\mae(X_t,X_0)=\frac{1}{\tau}\int_{-\tau}^{0}|X_t(s)-X_0(s)|ds$.
Bright colors in \cref{fig:heatmap} indicate large distances,
corresponding to large amplitude responses.  We notice an obvious
shift in behavior between $u=0.08$ and $u=0.09$ where the model goes
from exhibiting the small-amplitude periodic orbit to following the
large-amplitude solution.  This lower threshold is similar to the
observations when applying non-periodic insolation forcing (compare to
\cref{fig:heatmap_tau145}).

\subsection{Bifurcation analysis}
In order to examine the cause of the shift in behavior observed in \cref{fig:heatmap}, we first consider a numerical bifurcation analysis of DDE~\eqref{eq:periodic}.  The forcing period is kept constant at $T=4.1$ (corresponding to $41$ kyr).  We consider forcing amplitudes $u\in[0,0.75]$ and delays $\tau$ in the bistable region $\tau\in r_\mathrm{el}\approx[1.295,1.625]$ of the unforced system ($u=0$). 

\begin{figure}
\centering
\begin{subfigure}{0.48\textwidth}
	\centering
        \includegraphics[width=\textwidth]{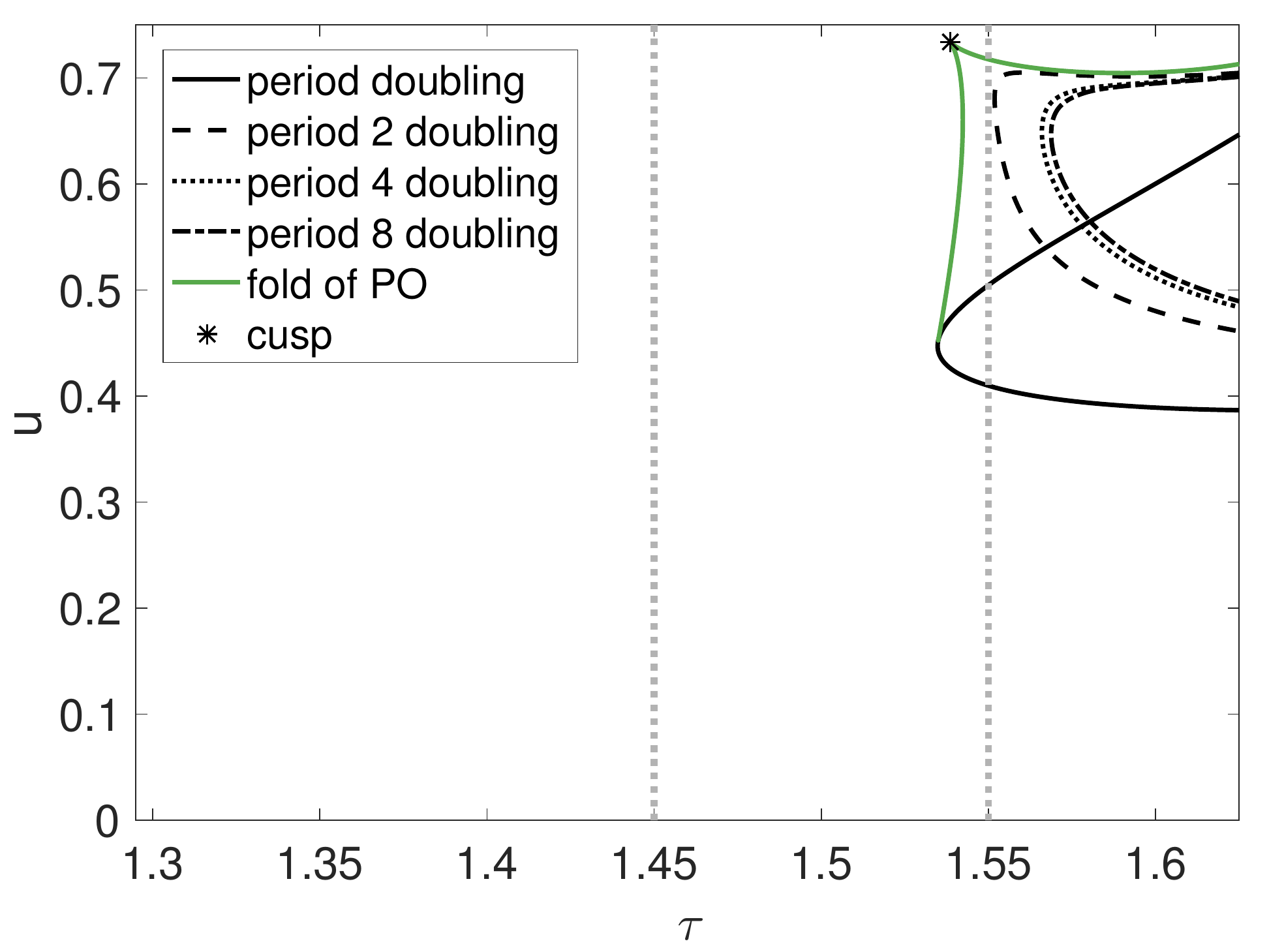}
        \caption{Bifurcation diagram for varying $\tau$ and $u$.}
        \label{fig:2parambif}
\end{subfigure}
\begin{subfigure}{0.48\textwidth}
        \centering
        \includegraphics[width=\textwidth]{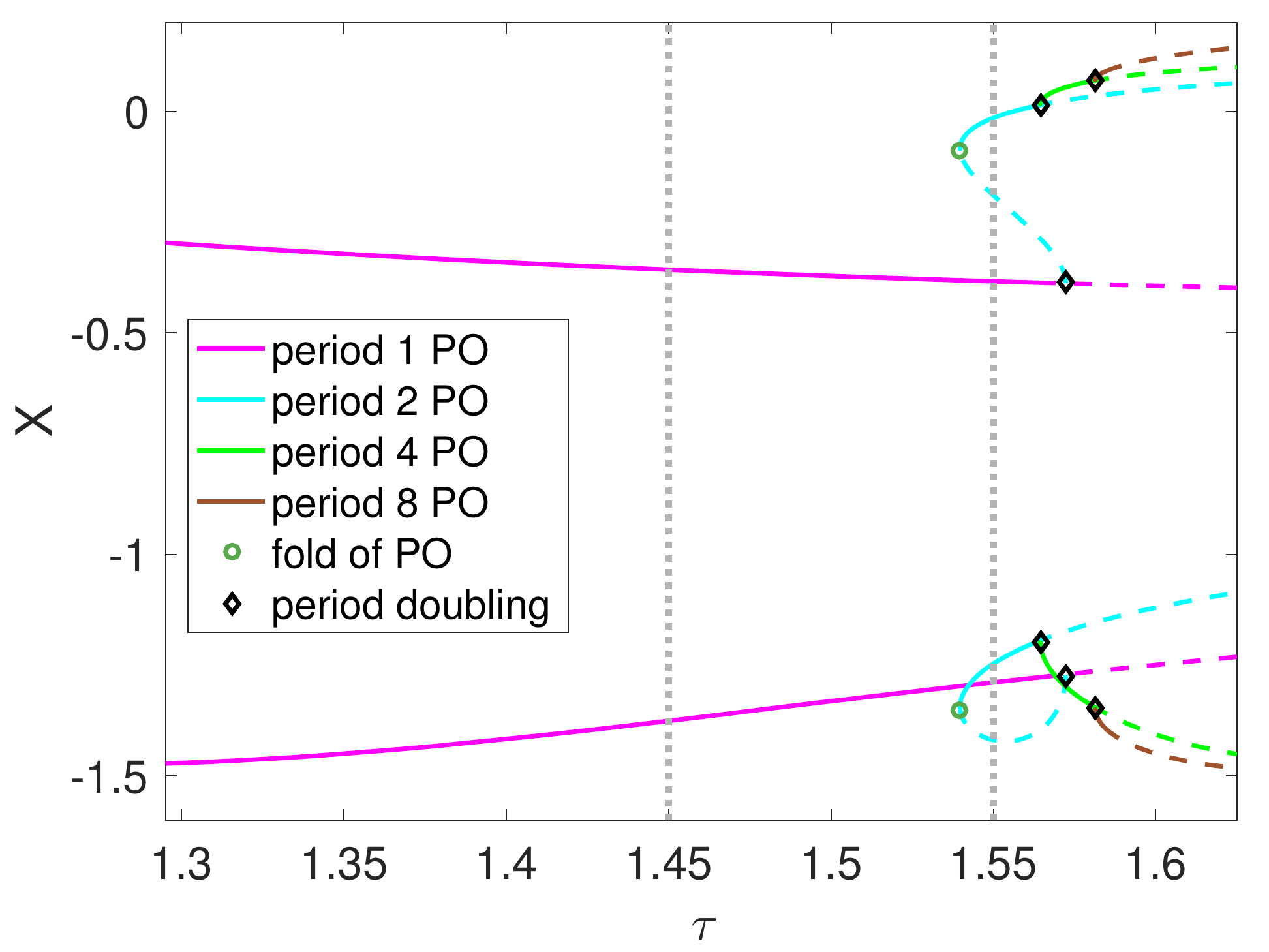}
        \caption{Cross section for $u=0.55$.}
        \label{fig:1parambif}
\end{subfigure}
\caption{Bifurcations of small-amplitude periodic orbit. Dotted vertical lines indicate values of $\tau$ used in \cref{sec:slow,sec:basin:planar,sec:phase}; other parameters:  $T=4.1$, $\tau=1.55$, $p=0.95$, $r=s=0.8$.}
\label{fig:periodic_bifdiag}
\end{figure}

\cref{fig:2parambif} shows the bifurcations of the small-amplitude periodic orbit. Bifurcations only occur for $\tau>1.53$ and $u>0.38$.  For a range of $u>0.4$ there exists a cascade of period doubling bifurcations for increasing $\tau$, evidence of which is also visible in~\cref{fig:heatmap}.  \cref{fig:1parambif} shows a cross section of the two-parameter bifurcation diagram \cref{fig:2parambif} along the horizontal line $u=0.55$ displaying the maximum and minimum of the periodic orbits on the $y$-axis.  We observe that the small-amplitude motion does \emph{not} experience any bifurcation for $u<0.3$. Moreover, because the bifurcations are restricted to large values of $\tau$, they cannot be used to explain the transition in time observed in \cref{fig:heatmap_tau145}, which is present for all $\tau$ throughout the bistable region; see \cite{quinn2017}.  Therefore periodic forcing, even with a slowly time-dependent modulated amplitude, is not sufficient to induce the MPT-like transition.

The large amplitude solution also goes through some bifurcations.  We do not show a detailed bifurcation analysis, but evidence of the collapse of the large amplitude solution can be seen in \cref{fig:heatmap}. These large responses are stable in a range of forcing amplitudes $u\in[0.09,0.15]$ in \cref{fig:heatmap}.  For $u\in[0.15,0.2]$ the trajectories make transient large-amplitude excursions before converging to a small-amplitude periodic orbit, which suggests a collapse of the (then chaotic) large-amplitude attractor.

As \cref{fig:2parambif} establishes, the observed transition in \cref{fig:heatmap} from
small- to large-amplitude oscillations at $u=0.09$ must have been caused by
some phenomenon other than a bifurcation. Since the unforced system is bistable for $\tau\in
r_\mathrm{el}$, we expect this bistability to persist for small
forcing amplitudes $u$. Thus, the initial history may cross from the
basin of attraction of the small-amplitude periodic orbit to the basin
of attraction of the large-amplitude response. Both example
trajectories in \cref{fig:ex_traj} started from the same initial
history but were computed with slightly different forcing amplitude
($u=0.08$ and $u=0.09$).  \cref{fig:heatmap} suggests that the constant initial history
$X_0:s\mapsto -0.5$ leaves the basin of attraction of the
small-amplitude periodic orbit at the lower threshold
$u\approx0.09$. 

\subsection{Dynamics on a two-dimensional slow manifold}
\label{sec:slow}
Since DDEs are infinite-dimension\-al, it is not feasible to determine
the basin of attraction in all dimensions.  However, previous studies
have proven results of Poincar{\'e}-Bendixson type (there exists a plane in $\R^2$ such that trajectories
cannot cross each other) for scalar
DDEs with monotone feedback \cite{mallet1990,mallet1996}. These are DDEs of the
form
\begin{align} \label{eq:mallet-DDE} &\dot{x}(t) = f(x(t),x(t-\tau)),
  \mbox{\ where\ $f(0,0)=0$ and $\delta y\,f(0,y)>0$ for
    all $y\neq0$, $\delta\in\{-1,1\}$,}
\end{align}
($\delta=1$ corresponds to positive delayed feedback, $\delta=-1$
corresonds to negative delayed feedback). 
The right-hand side
$f(X(t),X(t-\tau))$ in \eqref{eq:DDE} does not satisfy the feedback
conditions in \eqref{eq:mallet-DDE} since, for our right-hand side
$f$, $f(0,y)=-p y-s y^2$ changes sign also at $y=-p/s=-1.1875$, which
is reached in the unforced large-amplitude periodic orbit (see
\cref{fig:taubif_DDE}). However, the phase portrait in
\cref{fig:taubif_DDE} suggests that the unforced
DDE~\eqref{eq:periodic} (with $u=0$) has an attracting two-dimensional
slow manifold. We expect this manifold to persist for small forcing
amplitudes $u$. 

\begin{figure}
\centering
\begin{subfigure}{0.48\textwidth}
\includegraphics[width=\textwidth]{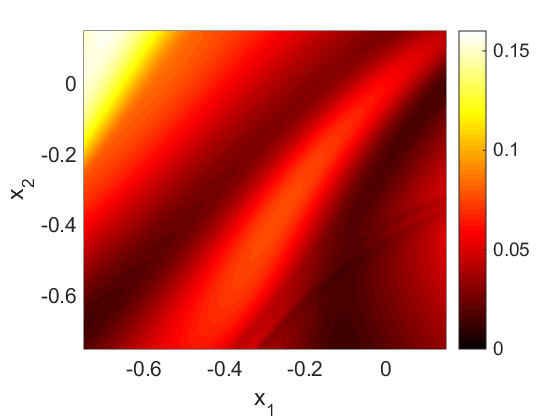}
\caption{\label{fig:spectralgap}}
\end{subfigure}
\begin{subfigure}{0.48\textwidth}
\includegraphics[width=\textwidth]{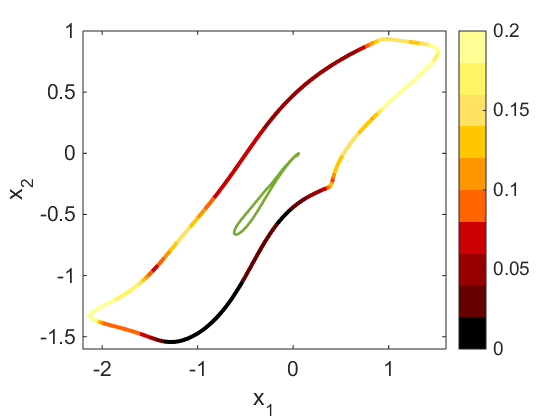}
\caption{\label{fig:singular}}
\end{subfigure}
\caption{Spectral gap of linearization of $M^1$ (a) on the domain $L\left([-0.75,0.15]\times[-0.75,0.15]\right)$ and (b) on the boundary for which the Jacobian of $RM^1L$ becomes singular. Color indicates the ratio between the third and second largest singluar values of $\partial M^1$.  The green boundary in figure (b) depicts the basin discussed in \cref{sec:basin:planar,sec:phase}. Parameters: $u=0.09$, $T=4.1$, $\tau=1.55$, $p=0.95$, $r=s=0.8$, $\phi=0$.}
\end{figure}

The apparent existence of an attracting two-dimensional slow
  manifold suggests that it may be possible to use an implicit
  computational dimension reduction introduced by Kevrekidis \emph{et
    al} as \emph{equation-free methods} (see reviews
  \cite{KS09,kevrekidis2010}).  The method can be used under the
  assumption that a high-dimensional system has a
    low-dimensional attracting slow manifold. The framework was
    originally developed for analysis of emergent macroscopic
    dynamics in stochastic or chaotic multi-particle
    simulations. Primarily, demonstrations of its use have focussed on
    analysis of equilibria or relative equilibria, e.g., bifurcation
    analysis (see \cite{thomas2016} for a recent general
    implementation) or control design \cite{SMK04}. In our case the
    high-dimensional system is the DDE \eqref{eq:periodic} (discretized by $N$ history
    points in practice) and the slow manifold is
    two-dimensional. Thus, the underlying problem is simpler than
    multi-particle systems as it has a well-understood time scale
    separation (demonstrated in \cref{fig:spectralgap}). However, we
    will construct the slow stable manifold of the period-one saddle
    periodic orbit, which is a more complex object than typically
    investigated with equation-free methods. The stable manifold will
    then pinpoint precisely the boundary of the basin of attraction
    inside the slow manifold.
  
  The idea of general equation-free framework is as follows.  One defines a \emph{lifting}
  function $L$ from the low-dimensional space $\R^2$ to the infinite-dimensional space
  of the DDE (e.g.\ $U=C([-\tau,0];\R)$).  One then evolves the dynamics on
  the infinite-dimensional space using a simulator of the high-dimensional system (the
  \emph{evolution} map $M$).  Here, this corresponds to solving DDE
  \eqref{eq:periodic} up to a time $t$ and extracting solution $X(s)$
  for $s\in[t-\tau,t]$.  Finally, one defines a \emph{restriction} function $R$
  to project the infinite-dimensional solution
  $X_t\in U=C([-\tau,0];\R)$ back into the low-dimensional space. In our case $RX_t$ is
  $(X_t(0),X_t(-\tau))=(X(t),X(t-\tau))\in \R^2$.  In summary, in our particular case, we have
\begin{align}\label{eq:liftdef}
  &\mbox{\emph{lifting}}& L&:\R^2\ni (x_1,x_2)\mapsto (X_0,\tilde X_0)\in \R\times \Lint^\infty([-\tau,0];\R)=:U^{\odot,*}\mbox{,}\\ &&&\nonumber 
  \mbox{\quad where $X_0=x_1$ and $\tilde X_0(s)=x_2$ for $s\in[-\tau,0]$,}\\
    &\mbox{\emph{evolution} map}&M_{t,t_0}&: U^{\odot,*}\ni (y_0,\tilde y)=X_{t_0}\mapsto X_{t_0+t}\in U \mbox{\quad for $t\geq\tau$,}\\
  \label{eq:restdef}
  &\mbox{\emph{restriction}}&R&:U\ni X \mapsto (X(0),X(-\tau))^T\in\R^2\mbox{.}
\end{align}
The range of $L$, called $U^{\odot,*}$, admits discontinuous bounded
segments and is a natural extension of the phase space
$U=C([-\tau,0];\R)$ of the DDE~\eqref{eq:periodic}. Trajectories
starting from $U^{\odot,*}$ return to the smaller phase space $U$
after time $\tau$.  See \cref{sec:ustar} for further comments.

The equation-free approach is then based on the two-dimensional map
\begin{equation} \label{eq:2Dmapping}
 RM_{t,0}L: \R^2 \ni (x_1,x_2)^T \mapsto (X(t),X(t-\tau))^T\in\R^2.
\end{equation}
The definition of the map in \eqref{eq:2Dmapping} means that for a
given pair $(x_1,x_2)\in\R^2$, we define the initial history of
\eqref{eq:periodic} as $X(0)=x_1$ and $X(s)=x_2$ for $s\in[-\tau,0)$,
simulate the DDE up to time $t\geq \tau$ from this history, and then use
$(X(t),X(t-\tau))^T$ as the result of the map. From this point onward we
will refer to $(x_1,x_2)$ as the argument of \eqref{eq:2Dmapping}.

In our computations we approximate elements $(X_0,\tilde X_0)\in
U^{\odot,*}$ by vectors $Y\in\R^N$, where $Y_k$ is an approximation of
 $\tilde X_0(-\tau(N-k)/(N-1))$ for $1\leq k<N$
and $Y_N=X_0$.  We use the discretized map $M$ based
on the Euler-Heun integration approximation with $h=0.01$, where a single step has
the form
\begin{align*}
  M_{h,t}&:\R^N\ni Y\mapsto \left(Y_2,\ldots,Y_N,Y_N+h(f_0+f_E)/2\right)^T\in\R^N\mbox{,}&&\mbox{where}\\
  f_0=&f(t,Y_N,Y_1,u)\mbox{,\quad} Y_N^0=Y_N+h f_0\mbox{,\quad}
  f_E=f(t+h,Y_N^0,Y_2,u)\mbox{,}&&\mbox{and}\\
  &f(t,x_1,x_2,u)=-p x_2+rx_1-sx_2^2-x_2^2x_1-uF_P(t)
\end{align*}
is the right-hand side of the DDE~\eqref{eq:periodic}. For larger time
spans we apply the composition rule
$M_{t+s,r}=M_{t,s+r}\circ M_{s,r}$ for $s,t\geq0$,
such that the discretization using $N-1=s/h$ steps converges to
  the continuous map $M_{t+s,s}$ uniformly for bounded
  $t\geq\tau$ and bounded intital values in $U^{\odot,*}$. We
restrict ourselves to stroboscopic maps $M_{t+s,s}$, where $t$ is
a multiple of the period: $t=kT$ with $k\in\Z$ and $T=4.1$, such that
we may write
\begin{displaymath}
  M^kY=M_{kT,0}(Y)
\end{displaymath}
for integers $k\geq0$.  The map $M^k:\rg L\to \dom R$ ($\rg L$ is the
range of $L$, $\dom R$ is the domain of definition of $R$) is autonomous and smooth, since
$M_{(k+j)T,jT}=M_{kT,0}$ for all integers $k\geq0$ and $j$,
and periodic forcing with period $T$. With this notation, $M^{k+j}$
equals $M^kM^j$. Compatible with the discretization of $M$, the
discretizations of lifting and restriction are
\begin{align*}
  L&:\R^2\ni(x_1,x_2)^T\mapsto
  (x_2,\ldots,x_2,x_1)^T\in\R^N\mbox{,}\\
  R&:\R^N\ni Y\mapsto(Y_N,Y_1)^T\in\R^2\mbox{.}
\end{align*}

Before proceeding with the application of the equation-free methods, we give numerical evidence that a two-dimensional slow invariant manifold is indeed present. \cref{fig:spectralgap} shows that the linearization of the map $M^1$
has a spectral gap after the first two eigenvalues such that $\partial
M^1$ is a small perturbation of a rank $2$ matrix for all $y$ in a
neighborhood of $L\left([-0.75,0.15]\times[-0.75,0.15]\right)$. This is numerical
evidence for the suspected time scale separation leading to a
two-dimensional slow manifold. We do not need to construct the slow
manifold explicitly, but rather may construct an approximate
two-dimensional map ${\cal M}_\ell$ from the slow manifold back to itself
implicitly, using coordinates in $\R^2$:
\begin{equation}
  \label{eq:mplanedef}
  {\cal M}_\ell:\R^2 \ni x\mapsto y\in\R^2 \mbox{,\quad where $y$ is the solution of\quad }
  RM^{\ell+1}Lx=RM^\ell Ly\mbox{.}
\end{equation}
The integer $\ell$ is the \emph{healing time} in the notation of
\cite{KS09,kevrekidis2010}. The map ${\cal M}_\ell$ approximates the
true stroboscopic map generated by the DDE~\eqref{eq:periodic} on the
slow manifold \cite{sieber2017convergence}. An intuitive
  explanation why the implicitly defined map $M_\ell$ is a valid
  approximation of the stroboscopic map on the slow manifold is given,
  for example, in \cite{Marschler2014b,sieber2017convergence}: the map $L$ maps $\R^2$ to a
  subspace that is assumed to be inside the basin of attraction of the
  attracting slow manifold. Thus, for both sides of the implicit
  definition \eqref{eq:mplanedef} the map $M^\ell L$ maps the element
  of $\R^2$ into the slow manifold (rather, very close to it, if $\ell$ is large enough). This
  map $M^\ell L$ is a diffeomorphism (a chart) between the slow
  manifold and $\R^2$. Calling $u_x=M^\ell Lx$ and $u_y=M^\ell Ly$,
  which are both objects inside the slow manifold in the
  high-dimensional space, the implicit definition \eqref{eq:mplanedef}
  requires $RMu_x=Ru_y$. If $R$ is a diffeomorphism between the slow
  manifold and $\R^2$ (a genericity condition) then this requirement
  implies that $Mu_x=u_y$. Thus, $M_\ell$, given by
  \eqref{eq:mplanedef}, approximates the map $M$ on the slow manifold
  in the coordinates given by the chart $M^\ell L$.

The approximation \eqref{eq:mplanedef} improves for increasing healing
time $\ell$ if lifting $L$ and restriction $R$ satisfy some genericity
conditions (implying that the map $R$ is a diffeomorphism between the
slow manifold and $\R^2$, and that $RM^\ell L:\R^2\to \R^2$ is a
diffeomorphism). The convergence result in
\cite{sieber2017convergence} does not require a large separation of
time scale, only a sufficiently large healing time. In our case
$\ell=1$ (a healing time of one period $T=4.1$) is sufficient: the
results only change by less than $10^{-2}$ when increasing $\ell$ to
$2$ (a large $\ell$ increases the condition number of
$\partial[RM^\ell L]$).

Within this persistent slow manifold the time-$T$ map
of the forced DDE~\eqref{eq:periodic} is a locally invertible
two-dimensional map. For two-dimensional maps the basin of attraction
for a periodic orbit is often bounded by the stable manifold of a
saddle periodic orbit. \cref{fig:spectralgap} justifies using a planar rectangle in
$\dom L$ to visualize the basins of attaction in the slow manifold.
Furthermore, since we can evaluate the stroboscopic map on the
two-dimensional slow manifold by using ${\cal M}_\ell$, we can employ
algorithms designed for the computation of stable manifolds of fixed
points in planar maps. In particular, we continue all three fixed
points present for zero forcing ($u=0$; see \cref{fig:taubif_DDE}) in
the parameter $u$, using the defining equation
\begin{displaymath}
  RM^\ell Lx_\mathrm{fix}=  RM^{\ell+1} Lx_\mathrm{fix}\mbox{,}
\end{displaymath}
which is a system of two equations for the two-dimensional variable
$x_\mathrm{fix}$ and the parameter $u$ (results will be shown for
healing time $\ell=1$). One of the fixed points is of saddle
type. The equation-free construction of a map from the slow
  manifold back to itself via \eqref{eq:mplanedef} permits us to
  extend specialized algorithms for planar maps such as the algorithm
  for the compuation of the stable manifold of a saddle fixed point
for maps that are not globally invertible, proposed by England
\emph{et al} \cite{england2004computing} and originally implemented
for two-dimensional maps in \texttt{DsTool}. Since the map
${\cal M}_\ell$ is implicitly defined, the algorithm as originally
implemented would require the solution of the nonlinear system
\eqref{eq:mplanedef} every time the map gets evaluated. This turns out
not to be necessary: we modify the stable manifold algorithm such that
it does not require any solution of a nonlinear system (see Appendix
\ref{sec:algo} for a brief explanation). The stable manifold of the
saddle fixed point will determine the basin of attraction for the
stable fixed point (the other fixed point is a source for all $u$) on
the slow manifold.

The large amplitude response (attracting initial conditions on the
other side of the stable manifold) is partially outside of the domain
of validity of the coordinates introduced by lifting $L$ and
restriction $R$. \cref{fig:singular} shows the curve in the
$\dom L$ plane along which the Jacobian of $RM^1L$ becomes singular,
which violates one of the assumptions made in the implicit definition
\eqref{eq:mplanedef} of ${\cal M}_1$.  Outside of this curve (where the large amplitude solution lies), our chosen plane is no longer valid.  However, the spectral gap values indicate that the dynamics may still be confined to a two-dimensional slow manifold.

\subsection{Basins of attraction and stable manifold in the plane}
\label{sec:basin:planar}
\begin{figure}
        \centering
        \begin{subfigure}{0.6\textwidth}
        \includegraphics[width=\textwidth]{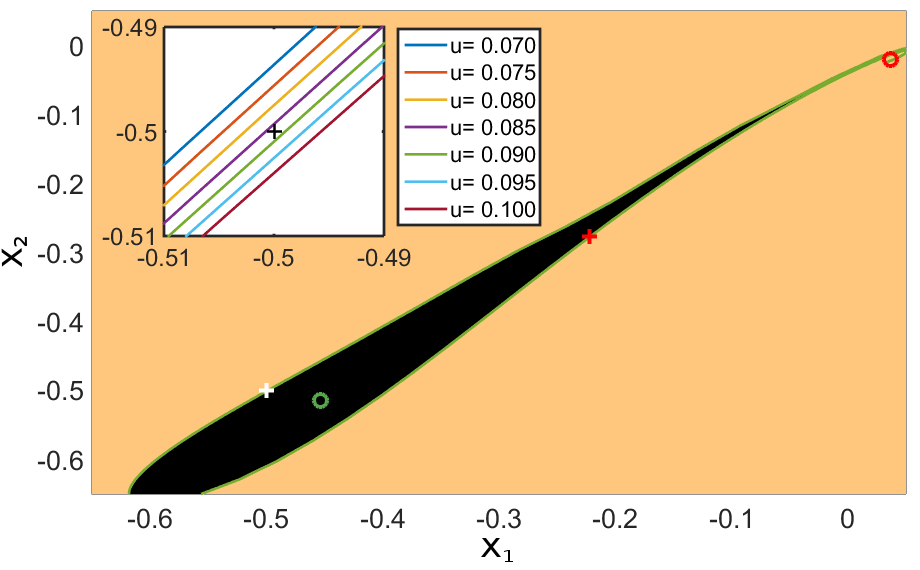}
	\caption{$\tau=1.55$}		
	\label{fig:basins}
		\end{subfigure}
		\begin{subfigure}{0.39\textwidth}
		\includegraphics[width=\textwidth]{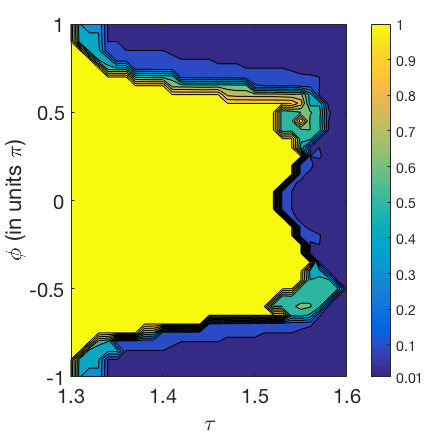}
	\caption{}
	\label{fig:phase_dep}		
		\end{subfigure}
\caption{\label{fig:basins:both} (a) Basin of attraction for the small amplitude stable periodic orbit for $\tau=1.55$.  Initial conditions in the black regions are attracted to the small amplitude stable periodic orbit (green circle) intersected with $\dom L$.  The white cross indicates the initial condition $(x_1,x_2)=(-0.5,-0.5)$ used to create \cref{fig:heatmap}.  The red circle is the unstable small amplitude periodic orbit, while the red cross represents the saddle periodic orbit. Figure zoom in top left shows stable manifold for different values of $u$ close to the initial condition $(x_1,x_2)=(-0.5,-0.5)$ (black cross). The basin was computed using simulations of \eqref{eq:periodic} with initial history $X(0)=x_1$ and $X(s)=x_2$ for $s\in[-\tau,0)$. The basin boundaries (colored lines) were computed using the method described in \cref{sec:algo}. (b) Threshold values for $u$ at which a transition to the
          large-amplitude response is observed as a function of delay
          $\tau$ and phase shift $\phi$ with initial history $X(s)=-0.5$ for $s\in[-\tau,0]$.  Values for which no transitions
          were observed are shown with $u=1$. Other parameters for both figures:
          $T=4.1$, $p=0.95$, $r=s=0.8$.} 
\end{figure}

Figure \ref{fig:basins} shows the basins in the
rectangle $[-0.65,0.05]\times[-0.65,0.05]$ for $u=0.09$.

The initial condition $x_1=x_2=-0.5$, corresponding to a 
constant initial history $X(s)=-0.5$ for $s\in[-\tau,0]$, was used in
the parameter scan for increasing forcing
amplitude $u$, depicted in \cref{fig:heatmap}.  This point is
indicated by a white cross in \cref{fig:basins}.  Black regions in
\cref{fig:basins} are initial histories $X(0)=x_1$ and $X(s)=x_2$ 
for $s\in[-\tau,0)$ that converge to the
stable small-amplitude periodic orbit in the center of the black
region. The beige region contains initial histories that escape to
the large amplitude response. The saddle fixed point is located near
$(-0.2,-0.3)$. Its stable manifold (in green) is the boundary
between the two basins of attraction.  The inset in \cref{fig:basins}
shows how the stable manifold of the saddle fixed point changes as the
forcing amplitude $u$ increases. In particular, we observe how the
initial condition $(-0.5,-0.5)^T$ is crossed by the stable manifold,
which shifts downward as $u$ increases.  An animation of the moving basin of attraction with increasing $u$ can be found in the supplementary material.

\subsection{Dependence on forcing phase}
\label{sec:phase}

As the basin of attraction in \cref{fig:basins} shows, the critical
forcing amplitude $u=0.09$ for the transition depends strongly on the
initial condition, which we chose as $X_0=-0.5$ ($=\mathrm{const}$)
for the heat map in \cref{fig:heatmap}. This dependence is
  specific to periodic forcing, since the transition for astronomical
  forcing shown in \cref{fig:heatmap_tau145} is not susceptible to
  small perturbations of the initial conditions (as pointed out in \cref{sec:astro}). Specifically, Quinn
  \emph{et al.}~\cite{quinn2017} reported negative finite-time
  Lyapunov exponents for times before $1000$ kyr BP such that a
  neighborhood of initial conditions leads to trajectories identical
  to those in \cref{fig:exbist}, starting from $X_0=-0.5$, after
  transients. As the periodically forced DDE possesses a slow
  manifold (see \cref{fig:spectralgap} for evidence), each initial
  condition in the infinite-dimensional phase space of the DDE leads
  to a trajectory that converges rapidly to a trajectory in the
  two-dimensional slow manifold. Thus, the initial condition for the
  periodically forced DDE should be in the plane shown in
  \cref{fig:basins} near $X_1=X_2=-0.5$ (since the astronomic forcing
  attracts us to this point and the slow manifold is attracting).  Therefore
  the only open question is
the
phase of the periodic forcing. We adjust the forcing equation
accordingly,
\begin{equation} \label{eq:forcing_periodic_phi}
F_P(t) = \sin((2\pi/T)\, t - \phi), \quad \phi \in [-\pi,\pi].
\end{equation}
The variable $\phi$ represents the phase shift of the forcing.  Note
that the bifurcation diagram \cref{fig:2parambif} is independent of
the forcing phase $\phi$.  However, \cref{fig:phase_dep} shows that
the phase affects the threshold value for the forcing amplitude $u$ at
which a transition to the large-amplitude response occurs with initial
history $X(s)=-0.5$ for $s\in [-\tau,0]$.
\cref{fig:phase_dep} shows contours of the smallest value of $u$ for
which we observe a transition to large-amplitude response in
simulations for different forcing phases $\phi$ and delays $\tau$ in
the bistable region $r_\mathrm{el}$ of \cref{fig:taubif_DDE}. For all
points in \ref{fig:phase_dep} we chose the initial history
  $X(s)=-0.5$ for $s\in [-\tau,0]$, however we note that changing the
  phase of the forcing is equivalent to considering different
  initial histories along a closed curve in the plane of
  \cref{fig:basins} through $(x_1,x_2)=(-0.5,-0.5)$. For some
parameter combinations the response is always small amplitude. In
these points we set the contour level to its maximum ($u=1$).  A
distinct boundary can be seen between parameter combinations that
exhibit transitions at low values of $u$ and those that do not.  For a
forcing phase $\phi = \pi$, a transition can always occur within the
bistable region.

In the supplementary material we show the effect of a phase shift on
the basin of attraction for the small-amplitude periodic orbit in the
plane $\dom L$.  This change of basin of attraction implies that, for
some phases $\phi$, trajectories from the initial history
$X(s)=-0.5$ for $s\in [-\tau,0]$ will converge to the small
  amplitude periodic orbit, while for other phases trajectories with
  the same initial history will converge to the large amplitude
  response.


\section{Conclusions}
\label{sec:conclusions}
We have presented a study of a scalar DDE model for the Pleistocene when subjected to periodic forcing.  The forcing period is 41 kyr, corresponding to the most prominent frequency seen in astronomical insolation variations.  With this particular forcing,  two stable solutions are observed at small forcing amplitudes: a small-amplitude periodic orbit and a large-amplitude, long-period response.  Both responses are a perturbation to the stable solutions of the autonomous system.  When the forcing amplitude is varied, we observe a threshold between convergence to the periodic orbit and convergence to the large-amplitude response.  A similar threshold is also seen with astronomic forcing presented in \cite{quinn2017}, where above a critical value of forcing strength $u$, transitions to a large-amplitude response occur.  Numerical bifurcation analysis of the small-amplitude periodic orbit shows that this threshold is outside of the region where any bifurcations occur, therefore leading us to conclude that changes in the basin of attraction for the periodic orbit and large amplitude response are the underlying mechanism.

In order to consider basins of attraction in the
  infinite-dimensional setting of DDEs, we give numerical
  evidence that there exists a two-dimensional slow manifold on which
a stable saddle, and unstable periodic orbit persist.  This
  justifies using equation-free methods to construct a stroboscopic
map on the slow manifold and track the locations of its fixed points.
We then apply a modified algorithm for computation of stable
  manifolds of saddle fixed points in planar maps to the slow manifold
  of the infinite-dimensional DDE without ever computing this slow
  manifold (original algorithm proposed by England \emph{et
    al.}~\cite{england2004computing}). Details and didactic
  implementation of the modified algorithm for implicitly defined maps
  are given in the supplementary material. This stable manifold forms
  the intersection between slow manifold and the boundary for the
basin of attraction of the stable periodic orbit.  We observe that the
stable manifold (and, with it, the basin boundary) moves for
increasing forcing amplitude.  It crosses our chosen initial history
for an amplitude around $u=0.09$.  We also scan the dependence of
  the basin of attraction on the phase of the forcing, which
  is equivalent to a change in the initial history function inside the
  slow manifold.  We conclude that the threshold behavior observed
for the periodically forced model can be attributed to a shifting
basin of attraction.

\paragraph{Connection to astronomic forcing --- a simplified scenario}
The results shown in
  \cref{fig:basins:both} provide a possible mechanism behind the
  observations for the astronomical forcing,
  recalled in \cref{sec:astro} and in \cref{fig:exbist}. In particular, the geometry shown in the inset in \cref{fig:basins} induces a transition similar to the one shown in \cref{fig:heatmap_tau145} for a simple step-wise amplitude modulation of the periodic forcing. If we replace the constant amplitude $u$ of the periodic forcing with a step-wise function of time,
  \begin{equation}\label{eq:umodulation}
    u(t)=
    \begin{cases}
      \ll1 &\mbox{for $t<750$kyr BP,}\\
      u_\mathrm{end} &\mbox{for $t\geq750$kyr BP,}
    \end{cases}
  \end{equation}
  then the response $X$ will be attracted to the area near the black
  cross in the inset of \cref{fig:basins} during the time
  $t\in[2000,750]$ kyr BP, independent of small perturbations of the
  initial condition (recall that the equilibrium for $u=0$ is at this
  point $X=-0.5$). At $t=750$ kyr BP, after the shift in $u$, the
  state $X$ will be outside of the basin of attraction of the fixed
  point for the new parameter value $u_\mathrm{end}$ (above the colored
  curve for $u$ in \cref{fig:basins}) if $u_\mathrm{end}$ is above a
  critical value ($\approx 0.09$). Thus, $X$ will escape toward the
  large amplitude response after $u=750$ kyr BP. Hence, for a
  step-wise periodic forcing with \eqref{eq:umodulation} one would
  make observations very similar to those listed in \cref{sec:astro},
  including the sudden transition with respect to time, with respect
  to $u_\mathrm{end}$, and the independence from small perturbations to
  the initial conditions. Scenarios where the shift of a basin of attraction 
  due to a parameter change causes a transition have been studied
  using the notion of \emph{rate-induced tipping}
  \cite{wieczorek2011excitability,ashwin2017parameter,ritchie2017probability,alkhayuon2018}. The
  scenario \eqref{eq:umodulation} corresponds to an infinite rate of
  change in the parameter (the forcing amplitude $u$) at time $t=750$
  kyr BP. For this scenario Ashwin \emph{et
    al.}~\cite{ashwin2017parameter} developed the concept of breakdown
  of basin forward stability, which generalizes the scenario described
  above for forcing \eqref{eq:umodulation}.

  The step-wise amplitude modulation \eqref{eq:umodulation} is simpler
  than the astronomic forcing (\cref{fig:forcing}, bottom
  panel). However, Quinn \emph{et al.}  \cite{quinn2017} showed for
  astronomic forcing that finite-time Lyapunov exponents (FTLEs) along
  the trajectories in \cref{fig:heatmap_tau145} remain negative for
  $u<0.22$ prior to transition time ($\sim750$ kyr BP). This implies
  that all trajectories from an open neighborhood of the unforced
  stable equilibrium are attracted to the same small-amplitude
  response prior to transition time.  Therefore the shift of basins of
  attraction is likely to be involved in the mechanism behind the
  observations in \cref{fig:exbist} and by Quinn \emph{et al.}
  \cite{quinn2017}.  

The numerical bifurcation analysis of the small-aplitude
  periodic orbit in \cref{fig:periodic_bifdiag} also illustrates that
  no bifurcation occurs for $\tau<1.53$.  This implies that even a
  gradual drift in the forcing amplitude (e.g. a slowly time
  dependent forcing amplitude $u(t)$ of the periodic forcing) is not
  enough to cause the temporal transition observed with the
  quasiperiodic forcing in \cref{fig:heatmap_tau145}. This supports
  the conclusions of \cite{quinn2017} that the phenomenon observed in
  \cref{fig:heatmap_tau145} is not related to a slow passage through a
  classical bifurcation, including bifurcations of periodic
  orbits. This is different from the idea presented in
  \cite{ditlevsen2018}, where in simple oscillator models under
  periodic forcing a generic transcritical bifurcation 
  is necessary for the system to transition between the
  smaller-amplitude 41 kyr cycle to the large amplitude 100 kyr
  response.  We therefore conjecture that the transition is an effect
  of the near-quasiperiodic modulation of the basic periodic forcing
    studied in this paper, which is present in the astronomical
  forcing, resulting in a more complicated version of the
    rate-induced tipping phenomenon obtained for the step-wise
    amplitude modulation \eqref{eq:umodulation}. This may require
  application of the general quasiperiodic theory of Fuhrmann \emph{et
    al.}  \cite{fuhrmann2017non} or the pullback attractor framework
  outlined in Chekroun \emph{et al.} \cite{chekroun2018}.

\appendix

\section{Extension of the phase space $U=C([-\tau,0];\R)$}
\label{sec:ustar}

This appendix gives a brief explanation for the extension of
  the phase space of DDEs, from which we permit initial conditions of
  the stroboscopic map $M$ in
  \cref{sec:slow,sec:basin:planar,sec:phase}. The space
  $\R\times\Lint^\infty([-\tau,0];\R)$ (called $U^{\odot,*}$ in the
  terminology of \cite{DGLW95}) is a natural extension of the phase
  space $U=C([-\tau,0];\R)$ of the DDE~\eqref{eq:periodic}.  The
  trajectories of the DDE~\eqref{eq:periodic} starting from initial
  values in $U^{\odot,*}$ admit discontinuous essentially bounded
  initial history segments $\tilde X_0:[-\tau,0]\mapsto\R$ and have
  $X_0$ as the right-side limit for $t\searrow0$. For an element
  $(X_0,\tilde X_0)$ of $U^{\odot,*}$, $X_0$ is usually called the
  \emph{head point}, while $\tilde X_0$ is the \emph{history segment}.
  The notation $\Lint^\infty([-\tau,0];\R)$ refers to the space of
  essentially bounded functions on $[-\tau,0]$ with essential maximum
  norm
  $\|\tilde X_0\|_0=\inf\left\{m\geq0:\leb\{t\in[-\tau,0]:|\tilde
    X_0(t)|\geq m\}=0\right\}$ ($\leb A$ is the Lebesgue measure of a
  set $A\subset\R$). The special role of the head point $X_0$ (compared to the
    remainder of the history segment $\tilde{X}_0$) becomes clear when re-stating a DDE
    using the equivalent variation-of-constants identity. In our special case of
    a DDE of the form $\dot X(t)=f(X(t),X(t-\tau),t)$ with fixed
    discrete delay $\tau$ this identity simplifies for $t\in[0,\tau]$
    to
    \begin{displaymath}
      X(t)=X_0+\int_0^tf(s,X(s),\tilde{X_0}(s-\tau))\mathrm{d}\,s\mbox{.}
    \end{displaymath}
    This identity makes clear that, for example, changing
    $\tilde{X}_0$ on a set of Lebesgue measure $0$ does not have any effect on
    the solution $X(t)$ for $t>0$, while changing $X_0$ does.

  As explained in the textbook \cite{DGLW95} trajectories starting
  from an element in the larger space $U^{\odot,*}$ return to the
  smaller phase space $U$ after time $\tau$ and the dependence
  $U^{\odot,*}\ni (X_0,\tilde{X_0}) \mapsto X_t\in U$ of the solution on its initial
  history is as regular as the right-hand side of the DDE
  \eqref{eq:periodic}. Thus,
\begin{displaymath}
  U^{\odot,*}\ni X\mapsto M(t,t_0;X)\in U
\end{displaymath}
is smooth for all $t\geq\tau$ and $t_0\in\R$.

\section{Summary of modification of stable manifold
  algorithm by England, Osinga and Krauskopf
  \cite{england2004computing}}
\label{sec:algo}
The supplementary material describes how one can modify the search circle (SC)
algorithm for stable manifolds in \cite{england2004computing} for maps
${\cal M}_\ell$ given implicitly through
\begin{equation}\label{app:imp:m}
    {\cal M}_\ell:\dom L\ni x\mapsto y\in\dom L\mbox{,\quad where $y$ is the solution of\quad }
  RM^{\ell+1}Lx=RM^\ell Ly\mbox{.}
\end{equation}
The supplement also contains a didactic implementation in
\texttt{Matlab} (\texttt{Gnu Octave} compatible) of the algorithm,
a demonstration script reproducing the manifold in
\cref{fig:basins}, and scripts reproducing some benchmark test
examples from \cite{england2004computing} (for the shear map, the
modified Ikeda map and the primary stable curve of the modified
Gumowski map; see \cite{england2004computing} for a review of the
properties and history of these examples).

The original SC algorithm \cite{england2004computing} grows the stable
curve of a map ${\cal M}_ \ell$ iteratively, approximating the stable
curve by a (linear) interpolation of a sequence of $k$ points
$S^k=(x^0,\ldots,x^k)$ in $\R^2$. At step $k+1$ a point $x^{k+1}$ is
added. The new point $x^{k+1}$ lies on a search circle arc with a
small adaptively chosen radius $\Delta$ around $x^k$. The point
$x^{k+1}$ is defined by the requirement that it is an intersection of
the image of this search circle under ${\cal M}_ \ell$ with the
previously computed curve $S^k$. In particular, the algorithm
\cite{england2004computing} does not rely on root-finding using Newton
iterations, but rather on a bisection to find the intersection between
image of the search circle and previous manifold. In principle, this
algorithm could be applied directly, if one solves the defining system
\eqref{app:imp:m}, $RM^{\ell+1}Lx=RM^\ell Ly$, for $y$ every time the
original algorithm applies its map (in our case ${\cal M}_\ell$) to a
point $x\in\R^2$.  However, a modification of the SC algorithm avoids
the need to solve the nonlinear equation \eqref{app:imp:m}. Instead of
a single sequence (and interpolating curve) $S^k$ one maintains two
curves, $S_L^k=(x_L^0,\ldots,x_L^k)$ in $\dom L$ and
$S_R^k=(x_R^0,\ldots,x_R^k)$ in $\rg R$, with $x_R^j=RM^\ell
Lx_L^j$. Then one searches for $x_L^{k+1}$ by finding an
intersection of the map $RM^{\ell+1}L$ image of the search circle of
radius $\Delta$ around $x_L^k$ in $\dom L$ with the curve $S_R^k$ in
$\rg R$ (note the power $\ell+1$ in the mapping). Once one found this
point $x_L^{k+1}$, for which $RM^{\ell+1}Lx_L^{k+1}\in S_R^k$, one
adds $x_L^{k+1}$ to $S_L^k$ and $x_R^{k+1}=RM^\ell Lx_L^{k+1}$ (note
the power $\ell$ of $M$ in the mapping) to $S_R^k$. Otherwise, the
same rules on angles and radii for acceptance of points apply as in
\cite{england2004computing}.

 \section*{Acknowledgments}
 We would like to thank the following people for their valuable discusisons and input: Timothy Lenton, Peter Ashwin, Peter Ditlevsen, Martin Rasmussen, Tobias J{\"a}ger, Flavia Remo, Hassan Alkhayuon, Paul Ritchie, and Damian Smug.

\bibliographystyle{siamplain}
\bibliography{references}
\end{document}